\documentclass[12pt,reqno]{amsart}
\allowdisplaybreaks[1]

\DeclareMathAlphabet{\mathpzc}{OT1}{pzc}{m}{it}

\author[]{Tomasz Przebinda\\ University of Oklahoma, Norman, OK, USA}
\title[The character and the wave front set correspondence]{The character and the wave front set\\ correspondence in the stable range}

\def\g{\mathfrak g}
\def\z{\mathfrak z}

\def\sp{\mathfrak {sp}}

\def\R{\mathbb{R}}
\def\C{\mathbb{C}}

\def\Ha{\mathbb{H}}

\def\z{\mathfrak z}

\def\ss1{\mathfrak s_{\overline 1}}

\def\hs1{\mathfrak h_{\overline 1}}

\def\Pg{\mathrm{P}}
\def\supp{\mathrm{supp}}

\def\Op{\mathrm{Op}}

\def\G{\mathrm{G}}
\def\N{\mathrm{N}}
\def\Zg{\mathrm{Z}}

\def\Qg{\mathrm{Q}}

\def\H{\mathrm{H}}
\def\M{\mathrm{M}}

\def\L{\mathrm{L}}
\def\Bbb{\mathbb}

\def\N{\mathrm{N}}

\def\H{\mathrm{H}}

\def\GL{\mathrm{GL}}

\def\Sp{\mathrm{Sp}}

\def\Og{\mathrm{O}}
\def\Ug{\mathrm{U}}

\def \t{\tilde}
\def \wt{\widetilde}

\def\V{\mathsf{V}}
\def\W{\mathsf{W}}

\def\V{\mathsf{V}}

\def\X{\mathsf{X}}
\def\Y{\mathsf{Y}}

\def\Ker{\mathrm{Ker}}

\def\End{\mathop{\hbox{\rm End}}\nolimits}

\def\det{\mathop{\hbox{\rm det}}\nolimits}

\def\Ad{\mathop{\hbox{\rm Ad}}\nolimits}
\def\Hom{\mathop{\hbox{\rm Hom}}\nolimits}

\def\tr{\mathop{\hbox{\rm tr}}\nolimits}

\def\sgn{\mathop{\hbox{\rm sgn}}\nolimits}

\def\lim{\mathop{\hbox{\rm lim}}\nolimits}

\def\supp{\mathop{\hbox{\rm supp}}\nolimits}

\def\Ss{\mathcal{S}}

\def\fonttitre{\textsf}
\newcounter{thh}
\newtheorem{thm}[thh]{\fonttitre{Theorem}}

\newtheorem{pro}[thh]{\fonttitre{Proposition}}
\newtheorem*{pro*}{\fonttitre{Proposition}}
\newtheorem{cor}[thh]{\fonttitre{Corollary}}
\newtheorem*{coro*}{\fonttitre{Corollary}}
\newtheorem{lem}[thh]{\fonttitre{Lemma}}

\newtheorem*{defi*}{\fonttitre{D??finition}}

\newtheorem*{nota*}{\fonttitre{Notation}}
\newenvironment{prf}{\begin{proof}}{\end{proof}}
\def\muet{ \ifthenelse{\equal{a}{b}}}
\def\nn{\nonumber}
\def\Z'{\Bbb{Z}'}

\def\biblio{\sloppy 
\bibliographystyle{alpha}            
\bibliography{article}}

\begin{document}



\maketitle

\begin{abstract}
We relate the distribution characters  and the wave front  sets of unitary representation for real reductive dual pairs of type I in the stable range.
\end{abstract}

\tableofcontents

\section{\bf Introduction. \rm}\label{Introduction}
In the late seventies Roger Howe formulated his theory of rank for irreducible unitary representations $\Pi$ of any connected cover  of the symplectic group $\Sp_{2n}(\R)$, see \cite{HoweRank}. The symplectic group has a maximal parabolic subgroup $\Pg$ with the Levi factor isomorphic to $\GL_n(\R)$ and the unipotent radical $\N$ isomorphic as a Lie group to the  space of the symmetric $n\times n$ matrices with the addition. In particular any connected cover of $\N$ splits.
The Spectral Theorem implies that the restriction of $\Pi$ to $\N$ is supported on the union of some $\GL_n(\R)$ - orbits in the dual of $\N$, which may be viewed as the space of the symmetric forms on $\R^n$. The rank of $\Pi$ is the maximal rank of a symmetric form in this support. 

A surprising result is that the representations $\Pi$ of rank $r<n$ are very special. The support of  $\Pi|_{\N}$ is a single $\GL_n(\R)$ - orbit of forms $\beta$ of signature $(p, q)$ with $p+q=r$. Furthermore, $\Pi$ factors through a double cover $\wt\Sp_{2n}(\R)$ of $\Sp_{2n}(\R)$
and remains irreducible when restricted to some other maximal parabolic subgroup $\wt\Pg_1\subseteq \wt\Sp_{2n}(\R)$. The Levi factor of $\Pg_1$ is isomorphic to $\GL_{r}(\R)\times\Sp_{2(n-r)}(\R)$ and the unipotent radical $\N_1$ is a two-step nilpotent group. The isometry group of a fixed form $\beta$ is isomorphic to $\Og_{p,q}\subseteq \GL_{r}(\R)$. According to \cite[Theorem 1.3]{Howesmall}, there is an irreducible unitary representation $\Pi'$ of $\wt\Og_{p,q}$ such that $\Pi|_{\wt\Pg_1}$ is induced from a representation involving $\Pi'$ of the subgroup $(\wt\Og_{p,q}\times\wt\Sp_{2(n-r)}(\R))\N_1\subseteq\wt \Pg_1$. 
The argument is based on the Stone von Neumann Theorem \cite{vonNeumannEindeutigkeit}, the theory of the Weil Representation \cite{WeilWeil} and the Mackey Imprimitivity Theorem, \cite{Mackey76}.

In particular the operators of $\Pi|_{\wt\Pg_1}$ are as explicit as the operators of $\Pi'$. However the remaining operators remain obscure. Fortunately there is a different description of the representations $\Pi$ and $\Pi'$.

The groups $(\Og_{p,q}, \Sp_{2(n-r)}(\R))$ form a dual pair in $\Sp_{2n}(\R)$ and there is Howe's correspondence for all real dual pairs $(\G, \G')$, \cite[Theorem 1]{HoweTrans}. As shown by Jian-Shu Li in his thesis, the representations $\Pi$ and $\Pi'$ are in Howe's correspondence. Li extended Howe's theory of rank to all dual pairs of type I and proved that it provides a bijection of representations of $\wt\G$ and $\wt\G'$ equal to Howe's correspondence, see \cite{Jian-ShuLiSingular} and \cite{Jian-ShuLiRank}. The condition of low rank is transformed to the dual pair being in the stable range, with $\G'$ - the smaller member. Now the operators $\Pi(g)$, $g\in \wt\G$, are much better understood because the Weil representation is known explicitly, see \cite{RangaRaoWeil} or section \ref{The  Weil representation} below, for a coordinate free approach.

Nevertheless an explicit description of all the $\Pi(g)$, $g\in \wt\G$, seems out of reach. Instead one may try to describe the distribution character $\Theta_\Pi$ of $\Pi$, \cite{HC-56}, in terms of  $\Theta_{\Pi'}$. This approach has a solid foundation, because for the dual pair $(\Ug_n, \Ug_n)$ the correspondence of the characters is governed by the Cauchy determinant identity, see \cite[Introduction]{PrzebindaCauchy}. In fact \cite[Definition 2.17]{PrzebindaCauchy} provides a candidate $\Theta'_{\Pi'}$ for $\Theta_\Pi$ in terms of $\Theta_{\Pi'}$. (For a more precise version see \cite[Formula (7)]{BerPrzeCHC_inv_eig}.) Let $\G'_1\subseteq \G'$ be the Zariski identity component.
Here is our first theorem. 
\begin{thm}\label{theteequalstheta'}
Suppose $(\G, \G')$ is a real irreducible dual pair of type I in the stable range with $\G'$ - the smaller member. Let $\Pi'$ be any genuine irreducible unitary representation of $\wt\G'$ and let $\Pi$ be the representation of $\wt\G$ corresponding to $\Pi'$. 
Assume that either $\G'=\G'_1$ or  $\G'\ne\G'_1$, but the restriction of $\Pi'$ to $\wt\G_1'$ is the direct sum of two inequivalent representations.
Then the restriction of $\Theta_\Pi$ to $\wt\G_1$ is equal to $\Theta'_{\Pi'}$. (For a character equality in the exceptional case see \eqref{theteequalstheta'0} below.)
\end{thm}
The proof looks as follows. As shown in \cite[Theorem 4]{BerPrzeCHC_inv_eig}, $\Theta'_{\Pi'}$ is an invariant eigendistribution. Hence, by Harish-Chandra Regularity Theorem, \cite[Theorem 2]{HC-63}, it suffices to know that the two distributions are equal on a Zariski open subset $\wt\G''\subseteq \wt\G$. This is verified using the method developed in  \cite{DaszPrzebindaInv} combined with a localization which requires the notion of a rapidly decreasing functions on $\wt\G$, as defined in \cite[7.1.2]{WallachI}.

Another invariant that tests our understanding of a representation is $WF(\Pi)$, the wave front set of $\Pi$. This notion, adapted from the theory differential operators, \cite[chapter 8]{Hormander},  was introduced to representation theory by Howe in \cite{HoweWave}. Since the wave front set of a representation of a reductive group is a union of nilpotent coadjoint orbits in the dual $\g^*$ of the Lie algebra $\g$ of $\G$, there are only finitely options for $WF(\Pi)$. Nevertheless it is surprisingly difficult to compute it. In part for that reason, Vogan introduced the notion of an associated variety of the representation (or rather of its Harish-Chandra module) in \cite{Vogan89}. As shown by Schmidt and Vilonen in \cite{SchmidVilonen2000}, the two notions are equivalent via the Sekiguchi correspondence of orbits, \cite{SekiguchiCorrespondence}. 

In order to state our second theorem, which expresses $WF(\Pi)$ in terms of $WF(\Pi')$, we need to recall that a dual pair $(\G, \G')$ is contained in the symplectic group $\Sp(\W)$, the isometry group of a nondegenerate symplectic form $\langle\cdot,\cdot\rangle$ on a finite dimensional vector space $\W$ over $\R$. Hence, there are moment maps $\tau_{\g}:\W\to\g^*$ and $\tau_{\g'}:\W\to\g'{}^*$ defined by
\begin{equation}\label{momentmaps}
\tau_{\g}(z)=\langle z(w),w\rangle \qquad (z\in\g,\ w\in\W)
\end{equation}
and similarly for $\g'$. 
%Also, recall \cite[Definition 8.1.2]{Hormander} that the wave front set of a distribution depends, up to the $\pm$ sign, on a choice of the Fourier transform. We make this %choice as in \cite[(1.14)]{PrzebindaCauchy}.
%%
\begin{thm}\label{theequality}
Suppose $(\G, \G')$ is a real irreducible dual pair of type I in the stable range with $\G'$ - the smaller member. Let $\Pi'$ be any genuine irreducible unitary representation of $\wt\G'$ and let $\Pi$ be the representation of $\wt\G$ corresponding to $\Pi'$. Then
\begin{equation}\label{theequality1}
WF(\Pi)=\tau_\g(\tau_{\g'}^{-1}(WF(\Pi')))\,.
\end{equation}
\end{thm}
We shall see in section \ref{Thewavefrontset} that Theorem \ref{theequality} follows from Theorem \ref{theteequalstheta'}, except when $\G'\ne \G'_1$ and the restriction of $\Pi'$ to $\wt\G'_1$ is irreducible. In that case, if $WF(\Pi')$ has more than one orbit of maximal dimension, we use a result of Loke and Ma, \cite{LockMaassocvar} combined with a theorem of Schmid and Vilonen,  \cite{SchmidVilonen2000}. In fact, \cite[Theorems A and D]{LockMaassocvar} prove the equality analogous to \eqref{theequality1} for all cases with the wave front set  replaced by the associated variety.
Therefore one is tempted to deduce \eqref{theequality1} from their result and from \cite{SchmidVilonen2000}. However, this is not straightforward, because Schmid and Vilonen \cite[page [1075]{SchmidVilonen2000} work with the groups that are the sets of the real points of a connected complex linear reductive groups. 
For a real ortho-symplectic dual pair either one member is a metaplectic group, which is not linear, or the other is an even orthogonal group, whose complexification is not connected. (Also, there are two Sekiguchi correspondences, see \cite{SekiguchiCorrespondence} and \cite[Proposition 6.6]{DaszKrasPrzebindaK-S2}, and the wave front set of a distribution depends, up to the $\pm$ sign, on a choice of the Fourier transform, see \cite[Definition 8.1.2]{Hormander}.)

One may probably circumvent \cite{LockMaassocvar} and  \cite{SchmidVilonen2000} by producing
the correct extension of $\Theta'_{\Pi'}$ from $\wt\G'_1$ to $\wt\G'$, but this would require a good understanding of the twisted orbital integrals, \cite{RenardTordues},  and is beyond the scope of this article.

The distribution $\Theta'_{\Pi'}$ is defined  also beyond  the stable range and does not depend on the unitarity of $\Pi'$. Furthermore, the Springer representations generated by the lowest terms in the asymptotic expansions of  $\Theta'_{\Pi'}$ and $\Theta_{\Pi'}$ behave as if \eqref{theequality1} were true beyond the stable range under some other mild assumptions, \cite[Theorem 1]{AubertKraskiewiczPrzebinda_real}. Therefore a generalization of the above two theorems  seems likely.
\section{\bf The  Weil representation. \rm}\label{The  Weil representation}
Fix a compatible positive complex structure $J$ on $\W$, i.e. $J\in \sp(\W)$ is such that $J^2=-1$, minus the identity in $\End(\W)$,  and the symmetric bilinear form $\langle J\cdot ,\cdot \rangle$ is positive definite. 
For an element $g\in \Sp(\W)$, let $J_g=J^{-1}(g-1)$. Then its adjoint with respect to the form $\langle J \cdot,\cdot\rangle$ is $J_g^*=Jg^{-1}(1-g)$. In particular  $J_g$ and $J_g^*$ have the same kernel. Hence the image of $J_g$ is
$
J_g\W=(\Ker J_g^*)^\perp=(\Ker J_g)^\perp
$, 
where $\perp$ denotes the orthogonal complement with respect to $\langle J \cdot,\cdot\rangle$.
Therefore, the restriction of $J_g$ to $J_g\W$ defines an invertible element. Thus it makes sense to consider $\det(J_g)_{J_g\W}^{-1}$, the reciprocal of the determinant of the restriction of $J_g$ to $J_g\W$. 
Let
%%
%\begin{equation}\label{metaplectic group}
\[
\wt{\Sp}(\W)=\{\t g=(g,\xi)\in \Sp(\W)\times \C,\ \ \xi^2=i^{\dim (g-1)\W}\det(J_g)_{J_g\W}^{-1}\}.
\]
%\end{equation}
%%
There exists a $2$-cocycle $C:\Sp(\W)\times \Sp(\W) \to \C$, so that $\wt{\Sp}(\W)$ is a group with respect to the multiplication 
%%
%\begin{equation}\label{multiplicationSp}
$(g_1,\xi_1)(g_2,\xi_2)=(g_1g_2,\xi_1\xi_2 C(g_1,g_2))$.
%\end{equation} 
%%
In fact, by \cite[Lemma 4.17]{AubertPrzebinda_omega},
\begin{equation}\label{abscocycle}
|C(g_1,g_2)|=\sqrt{\left|\frac{\det(J_{g_1})_{J_{g_1}\W}\det(J_{g_2})_{J_{g_2}\W}}{\det(J_{g_1g_2})_{J_{g_1g_2}\W}}\right|}
\end{equation}
and by \cite[Proposition 4.13 and formula (102)]{AubertPrzebinda_omega},
\begin{equation}\label{phasecocycle}
\frac{C(g_1,g_2)}{|C(g_1,g_2)|}=\chi(\frac{1}{8}\sgn(q_{g_1,g_2})),
\end{equation}
where $\chi(r)=e^{2\pi i r}$, $r\in \mathbb R$, is a fixed unitary character of the additive group $\R$ and  $\sgn(q_{g_1,g_2})$ is the signature of the symmetric form 
\begin{eqnarray*}
q_{g_1,g_2}(u',u'')=\frac{1}{2}\langle (g_1+1)(g_1-1)^{-1}u',u''\rangle&+&\frac{1}{2}\langle (g_2+1)(g_2-1)^{-1}u',u''\rangle\\ 
&&(u',u''\in (g_1-1)\W\cap (g_2-1)\W)\,.\nn
\end{eqnarray*}
By the signature of a (possibly degenerate) symmetric form we understand the difference between the maximal dimension of a subspace where the form is positive definite and the maximal dimension of a subspace where the form is negative definite. The group $\wt\Sp(\W)$ is known as the metaplecitc group.

Let $\W=\X\oplus\Y$ be a complete polarization. 
We normalize the Lebesgue measure on $\W$ and on each subspace of $\W$ so that the volume of the unit cube, with respect to the form $\langle J \cdot,\cdot\rangle$, is $1$. Since all positive complex structures are conjugate by elements of $\Sp(\W)$, this normalization does not depend on the particular choice of $J$.

Each tempered distribution $K\in \Ss^*(\X\times \X)$ defines an operator
$\Op(K)\in \Hom(\Ss(\X),\Ss^*(\X))$ by
%%
%\begin{equation}\label{Op}
\[
\Op(K)v(x)=\int_\X K(x,x')v(x')\,dx'.
\]
%\end{equation}
%%
Here $\Ss(\X)$ and $\Ss^*(\X)$ denote the Schwartz space on the real vector space $\X$ and the space of the tempered distributions on $\X$. 
The map $\Op: \Ss^*(\X\times \X)\to \Hom(\Ss(\X),\Ss^*(\X))$ is an isomorphism of linear topological spaces. This is known as the Schwartz Kernel Theorem, \cite[Theorem 5.2.1]{Hormander}. 

Fix the unitary character $\chi(r)=e^{2\pi i r}$, $r\in \mathbb R$, and
recall the Weyl transform 
\begin{eqnarray*}%\label{K}
&&\mathcal K:\Ss^*(\W)\to \Ss^*(\X\times \X)\,,\\
&&\mathcal K(f)(x,x')=\int_\Y f(x-x'+y)\chi\big(\frac{1}{2}\langle y, x+x'\rangle\big)\,dy \qquad (f \in \Ss(\W))\,.\nn
\end{eqnarray*}
Let
%\begin{equation} \label{eq:chicg}
\[
\chi_{c(g)}(u)=\chi\big(\frac{1}{4}\langle (g+1)(g-1)^{-1}u, u\rangle\big) \qquad (u=(g-1)w,\ w\in\W).
\]
%\end{equation}
In particular, if $g-1$ is invertible on $\W$, then 
$
\chi_{c(g)}(u)=\chi(\frac{1}{4}\langle c(g)u, u\rangle 
$ 
where $c(g)=(g+1)(g-1)^{-1}$ is the usual Cayley transform.
For $\t g=(g,\xi)\in\wt\Sp(\W)$ define
%%
%\begin{equation}\label{the omega}
\[
\Theta(\t g)=\xi,\qquad T(\t g)=\Theta(\t g)\chi_{c(g)}\mu_{(g-1)\W},\qquad \omega(\t g)=\Op\circ \mathcal K \circ T(\t g)\,,
\]
%\end{equation}
%%
where $\mu_{(g-1)\W}$ is the Lebesgue measure on the subspace $(g-1)\W$ normalized so that the volume of the unit cube with respect to the form $\langle J \cdot,\cdot\rangle$ is $1$. In these terms, $(\omega, \L^2(\X))$ is the Weil representation of $\wt\Sp(\W)$ attached to the character $\chi$.
In fact this is the Schr\"odinger model of $\omega$ attached to the complete polarization $\W=\X\oplus\Y$. Furthermore, $\Theta$ is the distribution character of $\omega$ and $T(\t g)$ is a normalized Gaussian. For future reference we set $\rho=\Op\circ\mathcal K$ and recall the following formula
\begin{equation}\label{tracegphi}
\tr \left(\omega(\t g)\rho(\phi)\right)=T(\t g)(\phi) \qquad (\t g\in \wt\Sp(\W), \phi\in \Ss(\W))\,.
\end{equation}
\section{\bf A mixed model  of the  Weil representation. \rm}\label{A mixed model of the  Weil representation}
In this section we recall the explicit  formulas for $\omega(\t g)$ for some particular elements $\t g$ of the metaplectic group.
For a subset $\M\subseteq \End(\W)$ let $\M^c=\{m\in \M\,:\,\det (m-1)\ne 0\}$ denote the domain of the Cayley transform in $\M$.
\begin{pro}\label{formula for M R}
Let $\M\subseteq \Sp(\W)$ be the subgroup of all the elements that preserve 
$\X$ and $\Y$. Set
\[
\det_\X^{-1/2}(\t m)=\Theta(\t m) |\det(\frac{1}{2}(c(m|_\X)+1))|^{-1} \qquad (\t m\in \wt\M^c).
\]
Then
%%
%\begin{equation}\label{formula for M 1 R}
\[
\left(\det_\X^{-1/2}(\t m)\right)^2=\det(m|_\X)^{-1} \qquad (\t m\in \wt\M^c)\,,
\]
%\end{equation}
%%
the function $\det_\X^{-1/2}\colon
\wt\M^c\to \C^\times$ extends to a continuous group homomorphism
\[
\det_\X^{-1/2}\colon\wt\M\to \C^\times
\]
and
%%
%\begin{equation}\label{formula for M 2 R}
\[
\omega(\t m)v(x)=\det_\X^{-1/2}(\t m) v(m^{-1}x) \qquad (\t m\in\wt\M,\ v\in \Ss(\X),\ x\in\X).
\]
%\end{equation}
%%
\end{pro}
Suppose $\W=\W_1\oplus\W_2$ is the direct orthogonal sum of two symplectic spaces. There are inclusions 
\begin{equation}\label{embedingofproductofgroups0}
\Sp(\W_1)\subseteq \Sp(\W),\ \ \ \Sp(\W_2)\subseteq \Sp(\W)
\end{equation}
defined by
\begin{eqnarray*}
g_1(w_1+w_2)&=&g_1w_1+w_2\,\\
g_2(w_1+w_2)&=&w_1+g_2w_2  \qquad (g_j\in \Sp(\W_j),\ w_j\in \W_j,\ j=1,2)\,.
\end{eqnarray*}
Furthermore, the map
\begin{equation}\label{embedingofproductofgroups}
\Sp(\W_1)\times \Sp(\W_2)\ni (g_1, g_2)\to g_1g_2\in \Sp(\W)
\end{equation}
is an injective group homomorphism.

Let us choose the compatible positive complex structure $J$ so that it preserves both $\W_1$ and $\W_2$. Then we have two metaplectic groups $\wt\Sp(\W_j)$, $j=1,2$. It is not difficult to see that the embeddings \eqref{embedingofproductofgroups0} lift to the embeddings
%%
%\begin{equation}\label{embedingofproductofgroups1}
\[
\wt\Sp(\W_1)\subseteq \wt\Sp(\W),\ \ \ \wt\Sp(\W_2)\subseteq \wt\Sp(\W).
\]
%\end{equation}
%%
Also, as is well known and easily follows from \eqref{abscocycle} and  \eqref{phasecocycle},
%%
%\begin{equation}\label{embedingofproductofgroups2}
\[
C(g_1, g_2)=1 \qquad (g_j\in \Sp(\W_j),\ w_j\in \W_j,\ j=1,2)\,.
\]
%\end{equation}
%%
Hence \eqref{embedingofproductofgroups} lifts to a group homomorphism
%%
%\begin{equation}\label{embedingofproductofgroups2}
\[
\wt\Sp(\W_1)\times \wt\Sp(\W_2)\ni (\t g_1, \t g_2)\to \t g_1\t g_2\in \wt\Sp(\W)\,,
\]
%\end{equation}
%%
with kernel equal to a two-element group. Moreover, in terms of the identification
%%
%\begin{equation}\label{embedingofproductofgroups3}
\[
\Ss(\W)=\Ss(\W_1)\otimes \Ss(\W_2)\,,
\]
%\end{equation}
%%
we have
%%
%\begin{equation}\label{embedingofproductofgroups4}
\[
T(\t g_1\t g_2)=T_{1}(\t g_1)\otimes T_{2}(\t g_2) \qquad (\t g_j\in \wt\Sp(\W_j),\ j=1,2)\,,
\]
%\end{equation}
%%
where $T_{j}(\t g_1)$ is the the normalized Gaussian for the space $\W_j$,  $j=1,2$. Hence,
%%
%\begin{equation}\label{embedingofproductofgroups5}
\[
\omega(\t g_1\t g_2)=\omega_{1}(\t g_1)\otimes \omega_{2}(\t g_2) \qquad (\t g_j\in \wt\Sp(\W_j),\ j=1,2)\,,
\]
%\end{equation}
%%
where $\omega_j$ is the Weil representation of $\wt\Sp(\W_j)$, $j=1,2$.

Suppose from now on that $\W_j=\X_j\oplus \Y_j$, $j=1,2$, are complete polarizations such that
\[
\X=\X_1\oplus\X_2\ \ \ \text{and}\ \ \ \Y=\Y_1\oplus \Y_2.
\]
Then, in particular, we have the following identifications
\begin{equation}\label{tensorproductidentification}
\Ss(\X)= \Ss(\X_1)\otimes  \Ss(\X_2)= \Ss(\X_1,\Ss(\X_2)).
\end{equation}
\begin{cor}\label{Maction}
Suppose $m\in \Sp(\W)$ preserves $\X_1$ and $\Y_1$. Denote by $m_1$ the restriction of $m$ to $\X_1$ and by $m_2$ the restriction of $m$ to $\Sp(\W_2)$.
Then for $v_1\in \Ss(\X_1)$, $v_2\in \Ss(\X_2)$, $x_1\in \X_1$ and $x_2\in \X_2$,
%%
%\begin{equation}\label{Maction1}
\[
\left(\omega(\wt {m_1}\wt {m_2})(v_1\otimes v_2)\right)(x_1+x_2)=\det_{X_1}^{-1/2}(\wt{m_1})v_1(m_1^{-1}x_1)(\omega_2(\wt{m_2}) v_2)(x_2)\,.
\]
Thus, in terms of \eqref{tensorproductidentification},
\[
\omega(\wt {m_1}\wt {m_2})v(x_1)=\det_{X_1}^{-1/2}(\wt{m_1})\omega_2(\wt{m_2})v(m_1^{-1}x_1) \qquad (v\in \Ss(\X_1,\Ss(\X_2)),\ x_1\in \X_1)\,.
\]
%\end{equation}
%%
\end{cor}
\begin{pro}\label{N1action}
%Suppose $n\in \Sp(\W)$ acts trivially on $\Y_1$ and on $\Y_1^\perp/\Y_1$. 
Suppose $n\in \Sp(\W)$ acts trivially on $\Y_1^\perp$. 
Then for
$v\in \Ss(\X_1,\Ss(\X_2))$ and $x_1\in \X_1$,
%%
%\begin{equation}\label{Naction1}
\[
%\omega(\t n)v(x_1)=\pm \chi_{c(-n)}(2x_1)\chi(\langle c(-n)x_1,x_2\rangle) v(x_1+x_2-c(-n)x_1)\,.
\omega(\t n)v(x_1)=\pm \chi_{c(-n)}(2x_1)v(x_1)\,.
\]
%\end{equation}
%%
%The $c(-n)x_1$ in $v(x_1+x_2-c(-n)x_1)$ stands for the $\X_2$-component of $c(-n)x_1$.
\end{pro}
\section{\bf The restriction of the Weil representation to the dual pair. \rm}\label{The restriction of the Weil representation to the dual pair}
The defining module $(\V, (\cdot,\cdot))$   for the group $\G$ is a finite dimensional left vector space $\V$ over a division algebra $\Bbb D=\R$, $\C$ or $\Ha$, with a possibly trivial involution, and a nondegenerate hermitian or skew-hermitian form $(\cdot,\cdot)$ such that $\G\subseteq \End_\Bbb D(\V)$ is the isometry group of that form. Similarly we have the defining module $(\V,' (\cdot,\cdot)')$ for the group $\G'$. The stable range assumption means that there is an isotropic subspace $\X_{(1)}\subseteq \V$ such that $\dim \V'\leq \dim X_{(1)}$. Select an isotropic subspace $\Y_{(1)}\subseteq \V$,  complementary to $X_{(1)}^\perp$, and let $\V_{(2)}\subseteq \V$ be the orthogonal complement of $\X_{(1)}\oplus \Y_{(1)}$, so that $\V=\X_{(1)}\oplus\V_{(2)}\oplus Y_{(1)}$.

The symplectic space  may be realized as $\W=\Hom(\V,\V')$ with 
\begin{equation}\label{explicitsymplecticform}
\langle w',w\rangle=\tr_{\Bbb D/\R}(w^*w'),
\end{equation}
where $w^*\in \Hom(\V',\V)$ is defined by
$(W,v')'=(v,w^*v')$, where $v\in \V$ and  $v'\in \V'$. The group $\G'$ acts on $\W$ by the post-multiplication and the group $\G$ by the pre-multiplication by the inverse.
Set $\X_1=\Hom(\X_{(1)},\V')$, $\Y_1=\Hom(\Y_{(1)},\V')$ and $\W_2=\Hom(\V_{(2)},\V')$. Then $\Y_1$ and $\X_1^\perp$ are complementary isotropic subspaces of $\W$ with respect to the symplectic form \eqref{explicitsymplecticform} and $\W_2$ is the orthogonal complement of $\W_1=\X_1+\Y_1$. We shall work in the mixed model of the Weil representation adapted to the decomposition $\W=\X_1\oplus \W_2\oplus \Y_1$, as explained in the previous section.

For any symmetric matrix $A\in \GL(\R^n)$ define
\[
\gamma(A)=  \frac{e^{\frac{\pi i}{4} \sgn(A) }}
{\sqrt{|\det A|}}\,.
\]
The real vector space $\Y_1$,  is  equipped with  the scalar product $\langle J\cdot,\cdot\rangle$.
Given $z\in\g$, the formula
$
q_z(y,y')=\frac{1}{2}\langle z y,y'\rangle
$
defines a symmetric bilinear form on $\Y_1$. Denote by $A_z$ the matrix of this form with respect to any orthonormal basis of $\Y_1$.
Denote by $i_{\Y_1}:\Y_1\to\X_1\oplus\W_2\oplus \Y_1$ the injection and by  $p_{\X_1}:\X_1\oplus\W_2\oplus \Y_1\to\X_1$ the projection. 
The matrix $A_z$ depends only on the map $p_{\X_1}zi_{\Y_1}:\Y_1\to \X_1$. The stable range assumption implies that we may choose $\X_{(1)}$ and $\Y_{(1)}$ so
the set of such elements $z$ in non-empty. We shall fix such a choice for the rest of this article and
let $\gamma(q_{p_{\X_1}zi_{\Y_1}})=\gamma(A_z)$.

The complete polarization $\W_1=\X_1\oplus \Y_1$ leads to the Weyl transform $\mathcal K_1:\Ss^*(\W_1)\to\Ss^*(\X_1\times\X_1)$. Hence $\mathcal K_1\otimes 1
:\Ss^*(\W)\to\Ss^*(\X_1\times\X_1\times \W_2)$. In order to shorten the notation we shall write $\mathcal K_1$ for $\mathcal K_1\otimes 1$. Explicitly
\[
\mathcal K_1(f)(x,x', w_2)=\int_{\Y_1} f(x-x'+y+w_2)\chi\big(\frac{1}{2}\langle y, x+x'\rangle\big)\,dy \qquad (f \in \Ss(\W), x,x'\in \X_1, w_2\in \W_2)\,.
\]
By computing a Fourier transform of a Gaussian, as in \cite[Theorem 7.6.1]{Hormander}, we obtain the following Lemma.
\begin{lem}\label{partialweyltransformfor g}
Let $z\in \g^c$ be such that $p_{\X_1}zi_{\Y_1}$ is invertible. 
Then for $x,x'\in\X_1$ and $w_2\in \W_2$ we have
\begin{eqnarray}\label{partialweyltransformfor g1}
\mathcal K_1(T(\wt{c(z)}))(x,x',w_2)&=&\Theta(\wt{c(z)})\gamma(q_{p_{\X_1}zi_{\Y_1}})\\
&&\chi_z(x-x')\chi_{(p_{\X_1}zi_{\Y_1})^{-1}}(x+x'-p_{\X_1}(z(x-x')+zw_2))\nn\\
&&\chi(\frac{1}{2}\langle zw_2,x-x'\rangle)\chi_z(w_2)\,.\nn
\end{eqnarray}
Let $h\in \G$ be the element that acts via multiplication by $-1$ on $\W_1$ and by the identity on $\W_2$. Suppose that in addition $\det(hc(z)-1)\ne 0$ and let $z_h=c(hc(z))$. Then
\[
\mathcal K_1(T(\wt{c(z_h)}))(x,x',w_2)=\det_{\X_1}^{-1/2}(\t h)\mathcal K_1(T(\wt{c(z)}))(x,-x',w_2)\,.
\]
(Here $\t h$ is one of the two elements in the preimage of $h$  chosen so that the right hand side is equal to the left hand side.)
\end{lem}
Here is a technical lemma, analogous to \cite[Lemma 4.3]{DaszPrzebindaInv}. Recall that for a test function $\Psi\in C_c^\infty(\wt\G)$
\[
T(\Psi)=\int_{\wt\G}\Psi(g)T(g)\,dg
\]
is a well defined tempered distribution on $\W$. Hence $\mathcal K_1(T(\Psi))$ is a tempered distribution on $\X_1\times \X_1\times \W_2$.
\begin{lem}\label{maintechicallemma}
Fix a euclidean norm $|\cdot|$ on the real vector space $\End(\V)$. There is a Zariski open subset $\G''\subseteq \G$ such that for $\Psi\in C_c^\infty(\wt\G'')$
the distribution $\mathcal K_1(T(\Psi))$ is a function on $\X_1\times \X_1\times \W_2$.
Moreover, for $N=0, 1, 2, ..$ there are constants $C_N$ such that for all $x,x'\in \X_1$ and all $w_2\in\W_2$,
\begin{eqnarray}\label{maintechicallemma1}
&&|\mathcal K_1(T(\Psi))(x,x',w_2)|\\
&\leq& C_N(1+|x^*x|+|x'{}^*x'|+|x^*x'|+|x'{}^*x|+|x^*w_2|+|x'{}^*w_2|+|w_2^*w_2|)^{-N}\,.\nn
\end{eqnarray}
\end{lem}
\begin{prf}
The function \eqref{partialweyltransformfor g1} is of the form $e^{i\frac{\pi}{2}\phi_{x,x',w_2}(z)}$, where
\begin{eqnarray*}
\phi_{x,x',w_2}(z)&=&\langle z(x-x'),x-x'\rangle\\
&+&\langle (p_{\X_1}z i_{\Y_1})^{-1}(x+x'-p_{\X_1}(z(x-x')+zw_2)), x+x'-p_{\X_1}(z(x-x')+zw_2)\rangle\\
&+&2\langle zw_2,x-x'\rangle+\langle zw_2,w_2\rangle\,.
\end{eqnarray*}
In order to simplify the computations we introduce the following notation
\begin{eqnarray*}
&&A=p_{\X_{(1)}}zi_{\X_{(1)}},\ \ \ B=p_{\X_{(1)}}zi_{\Y_{(1)}},\ \ \ C=p_{\Y_{(1)}}zi_{\X_{(1)}},\ \ \ F=C^{-1},\\ 
&&D=p_{\V_{(2)}}zi_{\Y_{(1)}},\ \ \ E=p_{\V_{(2)}}zi_{\X_{(1)}},\ \ \ z_2=p_{\V_{(2)}}zi_{\V_{(2)}}\,.
\end{eqnarray*}
By using the explicit description of the symplectic form,  \eqref{explicitsymplecticform}, and remembering that the Lie algebra $\g$ acts on $\W$ via minus the right multiplication, we can view the $A$, $B$, ..., $F$ as elements of $\End(\V)$, so that
\begin{eqnarray*}
-\phi_{x,x',w_2}(z)&=&\tr_{\Bbb D/\R}\Big( (x-x')^*(x-x')B\\
&+&(x+x'+(x-x')A+w_2E)^*(x+x'+(x-x')A+w_2E)F\\
&+&2(x-x')^*w_2D+w_2^*w_2z_2\Big)\,.
\end{eqnarray*}
The derivative of $-\phi_{x,x',w_2}(z)$ viewed as a function of the variables $A$, $B$, $F$, $D$, $E$, ${z_2}$ is given by
\begin{eqnarray*}
&-&\phi'_{x,x',w_2}(z)(\Delta_A, \Delta_B, \Delta_F, \Delta_D, \Delta_E, \Delta_{z_2})\\
&=&\tr_{\Bbb D/\R}\Big((x-x')^*(x-x')\Delta_B\\
&+&((x-x')\Delta_A)^*(x+x'+(x-x')A+w_2E)F\\
&+&(x+x'+(x-x')A+w_2E)^*(x-x')\Delta_AF\\
&+&(w_2\Delta_E)^*(x+x'+(x-x')A+w_2E)F+(x+x'+(x-x')A+w_2E)^*w_2\Delta_EF\\
&+&(x+x'+(x-x')A+w_2E)^*(x+x'+(x-x')A+w_2E)\Delta_F\\
&+&2(x-x')^*w_2\Delta_D+w_2^*w_2\Delta_{z_2}\Big)\,.
\end{eqnarray*}
Notice that $\Delta_AF=F(\Ad(F^{-1})\Delta_A)$. Also, by the structure of the Lie algebra $\g$, the variables $\Delta_A, \Delta_B, \Delta_F, \Delta_D, \Delta_E, \Delta_{z_2}$ are independent and fill out the corresponding vector spaces. The norm of the functional $\phi'_{x,x',w_2}(z)$ can be estimated from below by taking $\Delta_E=0$ and $\Delta_F=0$. Furthermore, all norms on a finite dimensional vector space are equivalent. 
Hence, with the appropriate choice of the norm $|\cdot |$ on $\End_\Bbb D(\V)$,
\begin{eqnarray}\label{maintechicallemma2}
|\phi'_{x,x',w_2}(z)|
&\geq&|(x-x')^*(x-x')|\\
&+&|(x-x')^*(x+x'+(x-x')A+w_2E)F|\nn\\
&+&|(x+x'+(x-x')A+w_2E)^*(x-x')F\Ad(F^{-1})|\nn\\
&+&2|(x-x')^*w_2|+|w_2^*w_2|\,.\nn
\end{eqnarray}
Using the operator norm inequality $|ab|\geq |a||b^{-1}|^{-1}$ and the fact that $|a^*|=|a|$ we see 
that
\begin{eqnarray*}
&&|(x-x')^*(x-x')|\geq |(x-x')^*(x-x')A||A|^{-1}\,,\\
&&|(x-x')^*(x+x'+(x-x')A+w_2E)F|\geq 
|(x-x')^*(x+x'+(x-x')A+w_2E)||F^{-1}|^{-1}\,,\\
&&|(x+x'+(x-x')A+w_2E)^*(x-x')F\Ad(F^{-1})|\\
&&\qquad \geq 
|(x-x')^*(x+x'+(x-x')A+w_2E)|\Ad(F)F^{-1}|^{-1}\,,\\
&&|(x-x')^*w_2|\geq |(x-x')^*w_2E||E|^{-1}\,.
\end{eqnarray*}

Hence,
\begin{eqnarray*}
&&|\phi'_{x,x',w_2}(z)|
\geq C(z)\big(|(x-x')^*(x-x')|+|(x-x')^*(x-x')A|+\\
&&\qquad 
|(x-x')^*(x+x'+(x-x')A+w_2E)|+|(x-x')^*w_2E|+|(x-x')^*w_2|+|w_2^*w_2|\big)\,,
\end{eqnarray*}

where
\[
C(z)=\min(\frac{1}{2}, \frac{1}{2}|A|^{-1}, |C|^{-1}+|\Ad(F)C|^{-1}, |E|^{-1})\,.
\]
Using the triangle 
inequality $|a|+|b|\geq |a\pm b|$ we see that
\begin{eqnarray*}
&&|(x-x')^*(x-x')A|+|(x-x')^*(x+x'+(x-x')A+w_2E)|+
|(x-x')^*w_2E|\\
&\geq&|-(x-x')^*(x-x')A+(x-x')^*(x+x'+(x-x')A+w_2E)-(x-x')^*w_2E|\\
&=&|(x-x')^*(x+x')|.
\end{eqnarray*}
So,
\begin{eqnarray*}
|\phi'_{x,x',w_2}(z)|\geq C(z)\big(|(x-x')^*(x-x')|+
|(x-x')^*(x+x')|+|(x-x')^*w_2|+|w_2^*w_2|\big)\,.
\end{eqnarray*}
All this is done under the condition on $z$ 
that the $C(z)$ is finite.  

Recall that Lemma \ref{partialweyltransformfor g} provides another expression for the function we would like to estimate, in terms $\phi'_{x,-x',w_2}(z_h)$. Our computation applied to $z_h$ shows that
\begin{eqnarray*}
|\phi'_{x,x',w_2}(z_h)|\geq C(z_h)(|(x-x')^*(x-x')|+|(x-x')^*(x+x')|+|(x-x')^*w_2|+|w_2^*w_2|)\,.
\end{eqnarray*}
Hence, by the triangle inequality again,
\begin{eqnarray*}
&&|\phi'_{x,x',w_2}(z)|+|\phi'_{x,-x',w_2}(z_h)|\\
&\geq&\min(C(z), C(z_h))\\
&&(|x^*x|+|x'{}^*x'|+|x^*x'|+|x'{}^*x|+|x^*w_2|+|x'{}^*w_2| +|w_2^*w_2|)\,.
\end{eqnarray*}
By the method of stationary phase (i.e. \cite[Theorem 7.7.1]{Hormander}) and Lemma \ref{partialweyltransformfor g}, the left hand side of \eqref{maintechicallemma1} is dominated by
\begin{eqnarray*}
&&\min((1+|\phi'_{x,x',w_2}(z)|)^{-N}, (1+|\phi'_{x,-x',w_2}(z_h)|)^{-N})\\
&\leq& (1+\frac{1}{2}(|\phi'_{x,x',w_2}(z)|+|\phi'_{x,-x',w_2}(z_h)|))^{-N}\\
&\leq& 2^N(\min(C(z), C(z_h))^{-N}\\
&&(|x^*x|+|x'{}^*x'|+|x^*x'|+|x'{}^*x|+|x^*w_2|+|x'{}^*w_2| +|w_2^*w_2||)^{-N}\,,
\end{eqnarray*}
which completes the proof, with $\G''$ equal to the image under the Cayley transform of the $z\in \g^c$ such that $z_h\in\g^c$ and both $C(z)$
and $C(z_h)$ are finite. 
\end{prf}
As an immediate consequence of Corollary \ref{Maction} and Proposition \ref{N1action} we deduce the following lemma.
\begin{lem}\label{G'actionandZ'action}
Let $\Zg\subseteq\G$ be the subgroup that acts trivially on $\Y_1^\perp$. Then for $\t n\in \wt{\Zg}$, $v\in \Ss(\X_1,\Ss(\X_2))$, $x_1\in \X_1$ and $\t g'\in\wt{\G'}$,
\begin{equation}\label{G'actionandZ'action1}
\omega(\t n)v(x_1)=\pm \chi_{c(-n)}(2x_1)v(x_1) \,,
\end{equation}
and
\begin{equation}\label{G'actionandZ'action2}
\omega(\t g')v(x_1)=\det_{X_1}^{-1/2}(\t g')\omega_2(\t g')v(g'{}^{-1}x_1) \,.
\end{equation}
\end{lem} 
\section{\bf The functions  $\Psi\in C_c^\infty(\wt\G'')$ act on $\mathcal H_\Pi$ via integral kernel operators. \rm}\label{via the integral kernel operators}
Given the polarization $\W_2=\X_2\oplus \Y_2$ we have the map
\[
\rho_2:\Ss^*(\W_2)\to \Hom(\Ss(\X_2), \Ss^*(\X_2))
\] 
as in \eqref{tracegphi}. Then 
\[
1\otimes \rho_2 : \Ss^*(\X_1\times\X_1\times \W_2)\to \Ss^*(\X_1\times\X_1)\otimes \Hom(\Ss(\X_2), \Ss^*(\X_2))\,.
\]
In order to shorten the notation we shall write  $\rho_2$ for $1\otimes \rho_2$ and $\mathcal K_1(T(\t g))(x,x')=\mathcal K_1(T(\t g))(x,x', \cdot)$. In these terms
\begin{equation}\label{complicatedomega}
\omega(\t g)v(x)=\int_{\X_1}\rho_2(\mathcal K_1(T(\t g))(x,x'))(v(x'))\,d\overset . x' \qquad (\t g\in\wt\G, v\in \Ss(\X_1,\Ss(\X_2))).
\end{equation}
Let $\X_1^{max}\subseteq \X_1$ be the subset of the surjective maps. The stable range assumption implies that this is a dense subset. 
For fixed $x, x'\in \X_1^{max}$ the operator norm of 
\begin{equation}\label{hsnorm}
\rho_2(\mathcal K_1(T(\Psi))(x,x')
\end{equation}
is bounded by the Hilbert-Schmidt norm, which is finite. Indeed, Lemma \ref{maintechicallemma} shows that
\[
\mathcal K_1(T(\Psi))(x,x', w_2)
\]
is a rapidly decreasing function of $x^*w_2$ and hence of $w_2$, because $x^*$, as a map from $\V'$ to $\X_{(1)}$ is injective. Therefore
\[
\mathcal K_1(T(\Psi))(x,x', \cdot)\in \L^2(\W_2)\,,
\]
which means that the Hilbert-Schmidt norm of \eqref{hsnorm} is finite.

In general, we denote by $\rho^c$ the representation contragredient to $\rho$ and by $\mathcal H_\rho$  a Hilbert space where $\rho$ is realized.

The group $\G'$ acts on $\X_1^{max}$, via the left multiplication, so that the quotient $\G'\backslash \X_1^{max}$ is a manifold. If $dx$ is a Lebesgue measure on $\X_1$, we shall denote by $d\overset . x$ the corresponding quotient measure on $\G'\backslash \X_1^{max}$.
Let $\mathcal H_\Pi$ be the Hilbert space of the functions $u:\G'\backslash \X_1^{max}\to \L^2(\X_2)\otimes \mathcal H_{\Pi'{}^c}$ such that for all $\t g'\in \wt\G'$
\begin{eqnarray}\label{realization of Pi}
u(g'x)=(\omega_2\otimes \det_{\X_1}^{-1/2}\Pi'{}^c)(\t g')u(x)\ \ \ \text{and}\ \ \ 
\int_{\G'\backslash \X_1^{max}}\parallel u(x)\parallel^2\,d\overset . x<\infty\,.
\end{eqnarray}
\begin{lem}\label{relizationasintegralkernels}
The representation $\Pi$ is realized on the Hilbert space $\mathcal H_\Pi$ and
for $\Psi\in C_c^\infty(\wt\G'')$, the operator $\Pi(\Psi)$ is given in terms of an integral kernel defined on $\X_1^{max}\times\X_1^{max}$ as follows
\[
(\Pi(\Psi)u)(x)=\int_{\G'\backslash \X_1^{max}}K_\Pi(\Psi)(x,x')u(x')\,d\overset . x' \quad(u\in \mathcal H_\Pi)\,,
\]
where
\begin{equation}\label{relizationasintegralkernels1}
K_\Pi(\Psi)(x,x')=\int_{\G'}\omega_2(\t g) \rho_2 (\mathcal K_1 (T(\Psi))(g^{-1}x,x',\cdot))\otimes \det_{\X_1}^{-1/2}(\t g)\Pi'{}^c(\t g)\,dg\,.
\end{equation}
Furthermore,
\begin{eqnarray}\label{relizationasintegralkernels2}
&&\tr K_\Pi(\Psi)(x,x')\\
&=&\int_{\G'} \int_{\W_2} T_2(\t g)(w_2) \mathcal K_1 (T(\Psi))(g^{-1}x,x', w_2)\det_{\X_1}^{-1/2}(\t g)\Theta_{\Pi'{}^c}(\t g)\,dw_2\,dg\,,\nn
\end{eqnarray}
where $\int_{\W_2} T_2(\t g)(w_2)\phi(w_2)\,dw_2$ stands for $T_2(\t g)(\phi)$.
\end{lem}
\begin{prf}
We proceed as in \cite[Proposition 4.8]{DaszPrzebindaInv}. Define a map 
\[
\Qg: \Ss(\X_1, \Ss(\X_2))\otimes \mathcal H_{\Pi'{}^c}\to \mathcal H_{\Pi}
\] 
by
\begin{eqnarray}\label{realization of Q}
\Qg(v\otimes \eta)(x)=\int_{\G'}(\omega\otimes\Pi'{}^c)(\t g)(v\otimes \eta)(x)\, dg\,.
\end{eqnarray}
Then \eqref{G'actionandZ'action2} shows that
\[
\Qg(v\otimes \eta)(x)=\int_{\G'}\omega_2(\t g)(v(g^{-1}x))\otimes\det_{\X_1}^{-1/2}(\t g)\Pi'{}^c(\t g)\eta\, dg\,.
\]
This last integral converges because $|g^{-1}x|$ is a constant multiple of the norm of $g$, as defined in \cite[2.A.2.4]{WallachI}. (The constant depends on $x$, which is fixed.) The argument used in the proof of Lemma 3.11 in \cite{DaszPrzebindaInv} shows that the range of $\Qg$ is dense in $\mathcal H_\Pi$. The action of  $\t g\in \wt\G$ on $\mathcal H_\Pi$ is defined via the the action of $\omega(\t g)$ on the $v$. Furthermore, with $\pi(\t g)=\det_{\X_1}^{-1/2}(\t g)\Pi'{}^c(\t g)$, we have
\begin{eqnarray*}
&&\Qg(\omega(\Psi)v\otimes\eta)(x)\\
&=&\int_{\G'}\omega_2(\t g)((\omega(\Psi)v)(g^{-1}x))\otimes \pi(\t g)\eta\,d g\\
&=&\int_{\G'}\int_{\X_1^{max}}\omega_2(\t g)\rho_2(\mathcal K_1(T(\Psi))(g^{-1}x,x'))(v(x'))\otimes \pi(\t g)\eta\,dx'\,d g\\
&=&\int_{\X_1^{max}}\int_{\G'}\omega_2(\t g)\rho_2(\mathcal K_1(T(\Psi))(g^{-1}x,x'))(v(x'))\otimes \pi(\t g)\eta\,d g\,dx'\\
&=&\int_{ \G'\backslash\X_1^{max}}\int_{\G'}\int_{\G'}\omega_2(\t g)\rho_2(\mathcal K_1(T(\Psi))(g^{-1}x,h^{-1}x'))(v(h^{-1}x'))\otimes \pi(\t g)\eta\,d g\,dh\,d\overset . x'\\
&=&\int_{ \G'\backslash\X_1^{max}}\int_{\G'}\int_{\G'}\omega_2(\t g\t h)\rho_2(\mathcal K_1(T(\Psi))(h^{-1}g^{-1}x,h^{-1}x'))(v(h^{-1}x'))\otimes \pi(\t g\t h)\eta\,d g\,dh\,d\overset . x'\\
&=&\int_{ \G'\backslash\X_1^{max}}\int_{\G'}\int_{\G'}\omega_2(\t g\t h)\omega_2(\t h)^{-1}\rho_2(\mathcal K_1(T(\Psi))(g^{-1}x,x'))(\omega_2(\t h)v(h^{-1}x'))\otimes \pi(\t g\t h)\eta\,d g\,dh\,d\overset . x'\\
&=&\int_{ \G'\backslash\X_1^{max}}\left(\int_{\G'}\omega_2(\t g)\rho_2(\mathcal K_1(T(\Psi))(g^{-1}x,x'))\otimes \pi(\t g)\,dg\right)\Qg(v\otimes\eta)(x')\,d\overset . x'\,,
\end{eqnarray*}
where by Lemma \ref{maintechicallemma} all the integrals are convergent.
This verifies \eqref{relizationasintegralkernels1}. 

The usual argument shows that $\mathcal K_1(T(\Psi))(g^{-1}x,x',w_2)$ is a differentiable function of $g$ and $w_2$ and that the derivatives are rapidly decreasing, as in Lemma \ref{maintechicallemma}. Hence \eqref{relizationasintegralkernels2} follows from \eqref{relizationasintegralkernels1} and \eqref{tracegphi}.
\end{prf}
\section{\bf The equality $\Theta_\Pi=\Theta'_{\Pi'}$. \rm}\label{The character}
Recall the group $\Zg$ defined in Lemma \ref{G'actionandZ'action}. For a Schwartz function $\psi\in\Ss(\z)$, on the Lie algebra $\z$ of $\Zg$, define a distribution $\psi_\Zg$ on $\wt\G$  by $\psi_\Zg=\t\psi\mu_\Zg$, where  $\t\psi(n)=\psi(c(-n))$, $n\in \Zg$, and $\mu_\Zg$ is the Haar measure on $\Zg$ viewed as a distribution on $\G$. Also, recall the space $\Ss(\G)$ of rapidly decreasing functions on $\G$, as defined in \cite[7.1.2]{WallachI}.
\begin{lem}\label{SanadS}
For any $\Psi\in C_c^\infty(\G)$ and any $\psi\in\Ss(\z)$, the convolution $\Psi*\psi_\Zg\in \Ss(\G)$.
\end{lem}
\begin{prf} 
Notice that for $z\in\z$, 
\[
-c(z)=(1+z)(1-z)^{-1}=(1+z)(1+z)=1+2z\,,
\]
because $z^2=0$. (Indeed, recall that  $z$ annihilates $\Y_1^\perp$. Since $(zx_1,y_1)=-(x_1,zy_1)=0$ we see that $z$ maps $\X_1$ to $\Y_1^\perp$, so $z^2=0$.)
Therefore we may choose the euclidean norm on the Lie algebra and the norm on the group, \cite[2.A.2.4]{WallachI}, so that $|c(z)|=|z|$, $z\in \z$.
Furthermore the map $\Zg\ni n\to c(-n)\in \z$ is a bijection with inverse $\z\ni z\to -c(z)\in \Zg$. 

Recall that
\[
\Psi*\psi_\Zg(a)=\int_\Zg\Psi(ab)\t\psi(b^{-1})\,d\mu_\Zg(b)\,.
\]
Let $C$ be a constant such that $|g|\leq C$ for all $g$ in the support of $\Psi$. Then
\begin{eqnarray*}
|\Psi*\psi_\Zg(a)|\leq \parallel\Psi\parallel_\infty\int_{|ab|\leq C}|\t\psi(b^{-1})|\,d\mu_\Zg(b)\leq \parallel\Psi\parallel_\infty 
 \int_{\frac{|a|}{C}\leq |b^{-1}|}|\t\psi(b^{-1})|\,d\mu_\Zg(b)\,.
\end{eqnarray*}
Since $\psi$ is rapidly decreasing,
\[
|\t\psi(b^{-1})|\leq C_N(1+|b^{-1}|)^{-N}\,.
\]
Furthermore,
\begin{eqnarray*}
\int_{\frac{|a|}{C}\leq |b^{-1}|}(1+|b^{-1}|)^{-N}\,d\mu_\Zg(b)\leq \int_{\frac{|a|}{C}\leq |b^{-1}|}(1+|b^{-1}|)^{-N/2}\,d\mu_\Zg(b) \left(1+\frac{|a|}{C}\right)^{-N/2}\,.
\end{eqnarray*}
Thus $\Psi*\psi_\Zg$ is rapidly decreasing. Further, we compute the left and right derivatives and get a similar estimate.
\end{prf}
Clearly, Lemma \ref{SanadS}, with the obvious modifications, holds for the groups $\wt\G$ and $\wt\Zg$, and we shall use it that way.
Define a Fourier transform
\[
\hat \psi(\zeta)=\int_\z \psi(z)e^{2\pi i\zeta(z)}\,dz\qquad (\zeta\in \z^*)
\]
and the moment map
\[
\tau_\z:\W\to \z^*,\ \ \ \tau_\z(w)(z)=\langle z w, w\rangle \qquad (w\in \W, z\in \z)\,.
\]
We shall see in \eqref{firstcharacterformula1} that the following lemma removes the ``deep'' stable range assumption from \cite{DaszPrzebindaInv}.
\begin{lem}\label{KPimodified}
For any $\Psi\in C_c^\infty(\wt\G'')$ and any $\psi\in\Ss(\z)$, 
\[
K_\Pi(\Psi*\psi_\Zg)(x,x')=2^{\dim \z} K_\Pi(\Psi)(x,x')  \hat \psi (\tau_\z(x')) \qquad (x,x'\in \X_1^{max})\,.
\]
\end{lem}
\begin{prf}
The formula \eqref{G'actionandZ'action1} implies that
\[
\mathcal K_1 (T(\Psi*\psi_\Zg))(x,x')=2^{\dim \z} \mathcal K_1 (T(\Psi))(x,x')\hat \psi (\tau_\z(x'))\,.
\]
Hence the lemma follows from \eqref{relizationasintegralkernels1}, because $\Zg$ commutes with $\G'$.
\end{prf}

In the remainder of this section we prove Theorem \ref{theteequalstheta'}.
Recall the distribution $\Theta'_{\Pi'}$ defined in \cite[Definition 2.17]{PrzebindaCauchy}. (For a more precise version see \cite[Formula (7)]{BerPrzeCHC_inv_eig}.) That invariant distribution was defined on smooth compactly supported functions, but that definition extends to $\Ss(\wt\G_1)$, without any modifications.

Let $\Psi\in C_c^\infty(\wt\G'')$ and let $\psi\in\Ss(\z)$, with $\supp\hat\psi$ compact. Denote by $\chi_{\Pi'}((-1\t ))$ the scalar by which $\Pi'((-1\t ))$ acts on the Hilbert space $\mathcal H_{\Pi'}$, so that $\Theta_{\Pi'}((-1\t )\t g)=\chi_{\Pi'}((-1\t ))\Theta_{\Pi'}(\t g)$. Also, recall \cite[sec. 4.5]{AubertPrzebinda_omega} the twisted convolution 
\[
\phi_1\natural\phi_2(w')=\int_\W\phi_1(w'-w)\phi_2(w)\chi(\frac{1}{2}\langle w,w'\rangle)\,dw \qquad (\phi_1, \phi_2\in \Ss(\W))\,,
\]
which extends by continuity to some tempered distributions so that, in particular,
\[
T(\t g_1)\natural T(\t g_2)= T(\t g_1\t g_2) \qquad (\t g_1, \t g_2\in \wt\Sp(\W))\,.
\]
Also, the formulas \eqref{relizationasintegralkernels1} and  \eqref{relizationasintegralkernels2} hold with $\Psi$ replaced by $\Psi*\psi_\Zg$, because $\wt\G'$ commutes with $\wt\G$.
Hence, with $\Phi=\Psi*\psi_\Zg$,
\begin{eqnarray}\label{firstcharacterformula1}
&&\Theta_\Pi(\Psi*\psi_\Zg)=\tr \Pi(\Phi)=\int_{\G'\backslash\X_1^{max}}\tr \mathcal K_\Pi (T(\Phi))(x,x)\,d\overset . x\\
&=&\int_{\G'\backslash\X_1^{max}}\tr \mathcal K_\Pi (T(\Phi))(-x,-x)\,d\overset . x\nn\\
&=&\int_{\G'\backslash\X_1^{max}}\int_{\G'}\int_{\W_2}T_2(\t g)(w_2) \mathcal K_1(T(\Phi))(-g^{-1}x,-x,w_2) \det_{\X_1}^{-1/2}(\t g)\Theta_{\Pi'}(\t g^{-1})\,dw_2\,dg\, d\overset . x\nn\\
&=&\chi_{\Pi'}((-1\t ))\int_{\G'\backslash\X_1^{max}}\int_{\G'}\int_{\W_2}T_2((-1\t )\t g)(w_2) \mathcal K_1(T(\Phi))(g^{-1}x,-x,w_2)\nn\\
&&\det_{\X_1}^{-1/2}((-1\t )\t g)\Theta_{\Pi'}(\t g^{-1})\,dw_2\,dg\, d\overset . x\nn\\
&=&\chi_{\Pi'}((-1\t ))\int_{\G'\backslash\X_1^{max}}\int_{\G'}\left(T_2((-1\t ))\natural T_2(\t g)\natural \mathcal K_1(T(\Phi))(g^{-1}x,-x,\cdot)\right)(0)\nn\\
&&\det_{\X_1}^{-1/2}((-1\t )\t g)\Theta_{\Pi'}(\t g^{-1})\,dg\, d\overset . x\nn\\
&=&\chi_{\Pi'}((-1\t ))\Theta_2((-1\t ))\int_{\G'\backslash\X_1^{max}}\int_{\G'}\int_{\W_2}\left(T_2(\t g)\natural \mathcal K_1(T(\Phi))(g^{-1}x,-x,\cdot)\right)(w_2)\nn\\
&&\det_{\X_1}^{-1/2}((-1\t )\t g)\Theta_{\Pi'}(\t g^{-1})\,dw_2\,dg\, d\overset . x\nn\\
&=&\chi_{\Pi'}((-1\t ))\Theta_2((-1\t ))\det_{\X_1}^{-1/2}((-1\t ))\int_{\G'\backslash\X_1^{max}}\int_{\G'}\int_{\W_2}\mathcal K_1(T_2(\t g)\natural T(\Phi))(g^{-1}x,-x, w_2)\nn\\
&&\det_{\X_1}^{-1/2}(\t g)\Theta_{\Pi'}(\t g^{-1})\,dw_2\,dg\, d\overset . x\nn\\
&=&\chi_{\Pi'}((-1\t ))\Theta_2((-1\t ))\det_{\X_1}^{-1/2}((-1\t ))\nn\\
&&\int_{\G'\backslash\X_1^{max}}\int_{\G'}\int_{\W_2}\int_{\Y_1} T(\t g)\natural T(\Phi)(x+x+y+w_2)\nn\\
&&\chi(\frac{1}{2}\langle y, x-x\rangle)\Theta_{\Pi'}(\t g^{-1})\,dy\,dw_2\,dg\, d\overset . x\nn\\
&=&\chi_{\Pi'}((-1\t ))\Theta((-1\t ))
\int_{\G'\backslash\X_1^{max}}\int_{\G'}\int_{\W_2}\int_{\Y_1} T(\t g)\natural T(\Phi)(x+y+w_2)\nn\\
&&\Theta_{\Pi'}(\t g^{-1})\,dy\,dw_2\,dg\, d\overset . x\nn\,,
\end{eqnarray}
where the functions under the integral are constant on the fibers of the covering map because we assume that $\Pi'$ is genuine,  and each consecutive integral is absolutely convergent, see Lemma \ref{KPimodified}. In fact, the integral over $(\G'\backslash\X_1^{max})\times \G'$ is also absolutely convergent. Also, $\Theta_{\Pi'}$ is indeed a function even if $\G'\ne \G'_1$, see \cite[Theorem 2.1.1]{Bouazizchar}.

Suppose first that $\G=\G'_1$.
Then we apply the Weyl integration formula for $\G'$
\[
\int_{\G'}f(g')\,dg'=\sum_{\H'}\frac{1}{|W(\H')|}\int_{\H'}\int_{\G'/\H'}g(g'h'g'{}^{-1})\,d\overset . g'\,|\Delta(h')|^2\,dh'
\]
to the integral over $\G'$ in \eqref{firstcharacterformula1} and see that
\begin{eqnarray}\label{firstcharacterformula2}
&&\Theta_\Pi(\Psi*\psi_\Zg)\\
&=&\chi_{\Pi'}((-1\t ))\Theta((-1\t ))
\sum_{\H'}\frac{1}{|W(\H')|}\int_{\H'\backslash\X_1^{max}}\int_{\H'}\int_{\W_2}\int_{\Y_1} T(\t h')\natural T(\Phi)(x+y+w_2)\nn\\
&&\Theta_{\Pi'}(\t h'{}^{-1})\,|\Delta(h')|^2\,dy\,dw_2\,dh'\, d\overset . x\nn\\
&=&\chi_{\Pi'}((-1\t ))\Theta((-1\t ))
\sum_{\H'}\frac{1}{|W(\H')|}\int_{\H'}\int_{\H'\backslash\X_1^{max}}\int_{\W_2}\int_{\Y_1} T(\t h')\natural T(\Phi)(x+y+w_2)\nn\\
&&\Theta_{\Pi'}(\t h'{}^{-1})\,|\Delta_{\G'}(h')|^2\,dy\,dw_2\,dh'\, d\overset . x\nn\\
&=&\chi_{\Pi'}((-1\t ))\Theta((-1\t ))
\sum_{\H'}\frac{1}{|W(\H')|}\int_{\H'}\int_{\H'\backslash\W^{max}}T(\t h')\natural T(\Phi)(w)\nn\\
&&\Theta_{\Pi'}(\t h'{}^{-1})\,|\Delta(h')|^2\,dh'\, d\overset . w\nn\\
&=&\chi_{\Pi'}((-1\t ))\Theta((-1\t ))
\sum_{\H'}\frac{1}{|W(\H')|}\int_{\H'}\Theta_{\Pi'}(\t h'{}^{-1})\,|\Delta(h')|^2\int_{\H'\backslash\W^{max}}\int_{\wt\G}\Phi(\t g)T(\t h'\t g)
\,d\t g\,dh'\, d\overset . w\nn\,.
\end{eqnarray}
Here we integrate over the regular  elements $\t h'\in\wt\H'$. For a fixed $\t h'$, the integral over $\H'\backslash\W^{max}$ is a distribution on the group $\wt\G$, which happens to be  the unique restriction of a distribution defined on the centralizer of the vector part of $\H'$ in $\wt\Sp(\W)$, as explained in \cite[Proposition 10]{PrzebindaCauchy}. Therefore \eqref{firstcharacterformula2} is equal to $\Theta'_{\Pi'}(\Psi*\psi_\Zg)$. However, as shown in \cite[Theorem 4]{BerPrzeCHC_inv_eig}, $\Theta'_{\Pi'}$ is an invariant eigendistribution.  (Being an eigendistribution is a local statement, so it does not depend on the class of the test functions.) Hence, by Harish-Chandra Regularity Theorem, \cite[Theorem 2]{HC-63}, we have the equality for a sufficient class of test functions to conclude that the two distributions on $\wt\G_1$ are equal.

Indeed, notice first that any invariant eigendistribution acts continuously on $\Ss(\wt\G)$ via an absolutely convergent integral. We see from the explicit formula for such a distribution restricted to a Cartan subgroup, \cite[Theorem 10.35]{knappLie2}, that this claim will be verified as soon as we check that the Harish-Chandra orbital integral corresponding to a Cartan subgroup $\wt\H\subseteq\wt\G$ maps $\Ss(\wt\G)$ continuously into the space $\Ss(\wt\H'')$ of the rapidly decreasing functions on the regular set $\wt\H''$ of $\wt\H$, as in \cite[Theorem 7.4.10(ii)]{WallachI}. If the $\wt\H$ is compact then the appropriate version of \cite[Lemma 7.4.5]{WallachI} carries over with an easier proof, because one does not have to invoke \cite[Lemma 7.4.3]{WallachI}. Hence, \cite[Lemma 7.A.4.2]{WallachI} implies that the last claim. Also, one sees from the inequalities \cite[(v), page 232]{WallachI} that the Harish-Chandra transform \cite[(v), page 231]{WallachI} maps $\Ss(\wt\G)$ continuously into the space of rapidly decreasing functions on the corresponding Levi factor. Hence, \cite[7.4.10(2), page 249]{WallachI} implies the claim for an arbitrary $\wt\H$. Next we notice that $\wt\G$ acts continuously on $\Ss(\wt\G)$ by translations, see the proof of Theorem 7.1.1 in \cite{WallachI}. We may choose the function $\psi_{\Zg}$ to be non-negative (squaring it results in the convolution on the Fourier side, which keeps the support compact) and have integral equal to $1$ and use the dilations on the Lie algebra $\z$ to construct the approximative identity, so that, with the notation of Lemma \ref{SanadS}, $\Psi*\psi_\Zg$ approaches $\Psi$ continuously in $\Ss(\wt\G)$. (This is a standard argument used for example in \cite[chapter 2, Theorem 2.1]{SteinIntroduction}.) Therefore our two invariant eigendistributions are equal on  $\Psi$. Thus they are equal on the Zariski open subset $\wt\G''\subseteq \wt\G$,. This suffices for the equality everywhere.

Now we consider the case $\G'_1\ne \G'$, where the Weyl integration formula is not valid.
The assumption that the restriction of $\Pi'$ to $\wt\G_1'$ is the direct sum of two inequivalent representations means that $\Theta_{\Pi'}$ is supported on $\wt\G_1'$. Each coset $\G'x\in \G'\backslash\X_1^{max}$ is the disjoint union of two disjoint $\G_1'$ cosets $\G'x=\G_1'x\cup(\G'\setminus\G_1')x$. The function we integrate in \eqref{firstcharacterformula1} has the same value on the two $\G'_1$ cosets but the integral over $\G'$ is equal to the integral over $\G'_1$, because $\Theta_{\Pi'}|_{\wt{\G'\setminus\G_1'}}=0$. Hence, modulo a factor of $2$, or renormalization of the measure on $\X_1$, the computation \eqref{firstcharacterformula1} goes through, with $\G'$ replaced by $\G'_1$. 

We continue studying the case $\G'_1\ne \G'$.
Suppose the restriction of $\Pi'$ to $\wt\G_1'$ is irreducible. Notice that $\wt\G'/\wt\G'_1$ is isomorphic to $\G'/\G'_1$. Hence the determinant may be viewed as a character
$\det:\wt\G'\to \C^\times$, trivial on $\wt\G_1'$. The representation $\Pi'\otimes\det$ is irreducible, is not equivalent to $\Pi'$ and has the same restriction to $\wt\G'_1$ as $\Pi'$. Let $\Pi_{\det}$ be the representation of $\wt\G$ corresponding to $\Pi'\otimes\det$. 
Since $\Theta_{\Pi'\otimes\det}=\Theta_{\Pi'}\det$, we see that the restriction of $\Theta_{\Pi'\oplus\Pi'\otimes\det}$ to $\wt{\G'\setminus\G_1'}$ is zero. Hence, the argument used in the proof of Theorem \ref{theteequalstheta'} shows that 
\begin{equation}\label{theteequalstheta'0}
\Theta_{\Pi\oplus\Pi_{\det}}=\Theta'_{\Pi'\oplus\Pi'\otimes\det}\,,
\end{equation}
where the right hand side is defined as before in terms or the Weyl integration formula for $\G'_1$.
\section{\bf The equality $WF(\Pi)=\tau_\g(\tau_{\g'}^{-1}(WF(\Pi'))$.\rm }\label{Thewavefrontset}
In this section we prove Theorem \ref{theequality}. If $\G'=\G'_1$, then the lowest term in the asymptotic expansion of $\Theta'_{\Pi'}$ is given in terms of the lowest in the asymptotic expansion of $\Theta_{\Pi'}$. This is immediate from \cite[Theorem 2.13]{PrzebindaCauchy}. Then \cite[Theorem 1.19]{PrzebindaCauchy} shows that by applying Fourier transform to both we get the desired orbit correspondence and \eqref{theequality1} follows from Theorem \ref{theteequalstheta'}. The same argument applies when $\G'\ne\G'_1$ and $\Pi'|_{\wt\G'_1}$ is the sum of two inequivalent representations. 

Suppose $\G'\ne\G'_1$ and the restriction of $\Pi'$ to $\wt\G'_1$ is irreducible. Then we have the equality \eqref{theteequalstheta'0} and the above argument shows that
\[
WF(\Pi\oplus\Pi_{\det})=\tau_\g(\tau_{\g'}^{-1}(WF(\Pi'\oplus\Pi'\otimes\det)))\,.
\]
Since the wave front set is computed at the identity and since the wave front set of the direct sum of representations is the union of their wave front sets \cite[Theorem 1.8 and Proposition 1.3(a)]{HoweWave}, we see that 
\[
WF(\Pi'\oplus\Pi'\otimes\det)=WF(\Pi')=WF(\Pi'\otimes\det)
\]
and
\[
WF(\Pi\oplus\Pi_{\det})=WF(\Pi)\cup WF(\Pi_{\det})\,.
\]
Thus
\begin{equation}\label{preliminaryequation}
WF(\Pi)\cup WF(\Pi_{\det})=\tau_\g(\tau_{\g'}^{-1}(WF(\Pi')))\,.
\end{equation}
Below we shall use \cite[Theorem 1.4]{SchmidVilonen2000} and \cite[Theorems A and D]{LockMaassocvar} in a minimal possible way. In particular we will not need these results if $WF(\Pi')$ is the closure of one orbit.

We know from \cite[Theorem 2.1.1]{Bouazizchar} that $\Theta_{\Pi'}$ has an asymptotic expansion near any semisimple point in $\wt\G'\setminus\wt\G'_1$. The corresponding asymptotic support at that point is contained in the wave front set at that point and hence, by \cite[Theorem 1.8]{HoweWave} in $WF(\Pi')$. Therefore the 
lowest possible homogeneity degree of the expansion at that point (a non-positive integer) is bounded below by the lowest possible homogeneity degree of the expansion at the identity. Therefore an obvious modification of \cite[Lemma 15(b)]{PrzebindaUnitary}, without the finite dimensionality assumption of the representation used there, holds and hence the argument of the proof of \cite[Theorem 7.8(b)]{PrzebindaUnitary} verifies the equality of the associate varieties of the primitive ideals
\[
Ass(I_{\Pi})=Ass(I_{\Pi_{\det}})\,.
\]
(After this improvement, \cite[Corollary E]{LockMaassocvar} is a particular case of \cite[Theorem 0.9]{PrzebindaUnitary}.)
Hence, by \cite[Theorem 4.1]{BarVogAs} and \cite{JosephAssociated}, the complexifications of $WF(\Pi)$ and $WF(\Pi_{\det})$ are equal.
By \cite[Theorem 8.1]{DaszKrasPrzebindaK-S2},  $\tau_\g(\tau_{\g'}^{-1}(WF(\Pi')))$ has the same number of nilpotent orbits of maximal dimension as $WF(\Pi')$.
If that number is $1$, then \eqref{preliminaryequation} shows that we can stop right here. Otherwise we rely on  \cite[Theorems A and D]{LockMaassocvar} and \cite[Theorem 1.4]{SchmidVilonen2000}, as explained below.

Since the restriction of $\Pi'$ to $\wt \G'_1$ is irreducible, one may take an irreducible subrepresentation of a maximal compact subgroup of $\wt \G_1$ in the definition of the good filtration leading to the associated variety of the Harish-Chandra module of $\Pi'$. Hence the associated varieties of the Harish-Chandra module of $\Pi'$ viewed as a representation of $\wt\G'$ or $\wt \G'_1$ are equal. The same argument applies to $\Pi'\otimes\det$. Thus we have the equality of the associated varieties $AV(\Pi')=AV(\Pi'\otimes\det)$. Then \cite[Theorems A and D]{LockMaassocvar} shows that $AV(\Pi)=AV(\Pi_{\det})$ . 

Since $\G'$ is an orthogonal group with the defining module of an even dimension, the covering $\wt\G\to \G$ is trivial, so \cite{SchmidVilonen2000} applies to $\wt\G$. Also, the group $\wt\G'_1$ is linear and the complexification of $\G'_1$ is connected. Hence \cite{SchmidVilonen2000} applies to $\wt\G'_1$. Therefore \cite[Theorem 1.4]{SchmidVilonen2000} together with \cite[Theorem 8.1]{DaszKrasPrzebindaK-S2} show that $WF(\Pi)=WF(\Pi_{\det})$, and we are done.

\biblio
\end{document}